\documentclass[centertags]{article}

\usepackage{amsmath}
\usepackage{amsthm}
\usepackage{amscd}
\usepackage{amssymb}
\usepackage{verbatim}
\title{Groups of homeomorphisms of one-manifolds, I: actions of 
nonlinear groups}
\author{Benson Farb\thanks{Supported in 
part by the NSF and by the Sloan Foundation.}\ \  and 
John Franks\thanks{Supported in part by NSF grant DMS9803346.}}
\date{June 20, 20001}
\newtheorem{theorem}{Theorem}[section]

\newtheorem{exmple}[theorem]{Example}

\newtheorem{thm}{Theorem}[section] 
\newtheorem{cor}[thm]{Corollary}
\newtheorem{lem}[thm]{Lemma} 
\newtheorem{prop}[thm]{Proposition}
\newtheorem{defn}[thm]{Definition}

\newcommand{\R}{{\bf R}}
\newcommand{\Z}{{\bf Z}}
\newcommand{\T}{S}

\DeclareMathOperator{\Fix}{Fix}
\DeclareMathOperator{\Int}{Int}
\DeclareMathOperator{\Per}{Per}
\DeclareMathOperator{\Homeo}{Homeo}

\DeclareMathOperator{\Diff}{Diff}

\DeclareMathOperator{\Aut}{Aut}
\DeclareMathOperator{\Out}{Out}
\DeclareMathOperator{\BS}{BS}
\DeclareMathOperator{\Mod}{Mod}
\DeclareMathOperator{\PSL}{PSL}
\DeclareMathOperator{\supp}{supp}

\begin{document}
\maketitle
\begin{abstract}
This self-contained paper is part of a series \cite{FF2,FF3} 
on actions by diffeomorphisms of infinite groups on compact manifolds.  
The two main results presented here are:
\begin{enumerate}
\item Any homomorphism of (almost any) 
mapping class group or automorphism group of a free group into
$\Diff_+^r(S^1), r\geq 2$ is trivial. For $r=0$ Nielsen showed that in 
many cases nontrivial (even faithful) representations exist.  Somewhat
weaker results are proven for finite index subgroups.

\item We construct a finitely-presented 
group of real-analytic diffeomorphisms of $\R$ which is not residually
finite. 
\end{enumerate}
\end{abstract}

\section{Introduction}

In this paper we consider infinite groups acting by
diffeomorphisms on one-dimensional manifolds.  For lattices in higher
rank semisimple Lie groups such actions are essentially completely
understood: 

\begin{theorem}[Ghys \cite{Gh}, Burger-Monod \cite{BM}]
\label{theorem:gbm}
Let $\Gamma$ be a lattice in a simple Lie group of $\R$-rank at least
two.  Then any $C^0$-action of $\Gamma$ on $S^1$ has a finite orbit, 
and any $C^1$-action of $\Gamma$ on $S^1$ must factor
through a finite group. 
\end{theorem}

Theorem \ref{theorem:gbm} is the solution in dimension one of Zimmer's
program of classifying actions of higher rank 
lattices in simple Lie groups on
compact manifolds.  In \cite{La} (see in particular \S{8}), Labourie
describes possible extensions of Zimmer's program to other ``big''
groups.  A. Navas \cite{Na} has recently proven that any $C^{1+\alpha},
\alpha>1/2$ action of a group with Kazhdan's property T factors through
a finite group. In this paper we consider three basic examples of nonlinear
groups: mapping class groups of surfaces, (outer) automorphism groups of
free groups, and Baumslag-Solitar groups.

\bigskip
\noindent
{\bf Mapping class groups and automorphism groups of free groups. }
Let $\Mod(g,k)$ denote the group of isotopy classes of diffeomorphisms
of the genus $g$ surface with $k$ punctures.  Unlike the case of 
lattices (Theorem \ref{theorem:gbm}), the group 
$\Mod(g,1), g\geq 1$ does have a faithful, $C^0$-action on $S^1$ without
a global fixed point.  This is a classical result of Nielsen.  
However, imposing a small amout of regularity changes the situation 
dramatically.

\begin{theorem}[Mapping class groups]
\label{theorem:mcgroups}
For $g\geq 3$ and $k=0,1$, any $C^2$ action of $\Mod(g,k)$ 
on $S^1$ or on $I=[0,1]$ is trivial.  
\end{theorem}

The $S^1$ case of Theorem \ref{theorem:mcgroups} was announced by E. Ghys 
several years ago (see \S{8} of \cite{La}), and for real-analytic
actions was proved by Farb-Shalen \cite{FS}.  

Another class of basic examples of ``big groups'' are automorphism
groups of free groups.  Let $\Aut(F_n)$ (resp. $\Out(F_n)$) denote the
group of automorphisms (resp. outer automorphisms) of the free group
$F_n$ of rank $n$. It is known that $\Aut(F_n),n>2$ is not linear
\cite{FP}.  

The techniques we develop to prove the results above allow us to prove
the following.

\begin{theorem}[Automorphism groups of free groups]
\label{theorem:autogroups}
For $n \ge 6$, any homomorphism from $\Aut( F_n)$ to
$\Diff_+^2(\T^1)$ factors through $\Z/2\Z$.  Any homomorphism from
$\Aut( F_n)$ to $\Diff_+^2(I)$ is trivial.  The analagous results
hold for $\Out(F_n)$.
\end{theorem}

M. Bridson and K. Vogtmann \cite{BV} have recently proven a 
stronger result: any homomorphism from $\Aut(F_n)$ to a group which
does not contain a symmetric group must have an image which is at most
$Z/2Z.$ However, their proof uses torsion elements in an essential way, 
hence does not extend to finite index subgroups.  
On the other hand our techniques have implications for finite
index subgroups of $\Aut(F_n)$ and $\Mod(g,k)$; such subgroups are 
typically torsion free.

\begin{theorem}[Finite index and other subgroups]
\label{theorem:finite:index}
Let $H$ be any group of $C^2$ diffeomorphisms of $I$ or $S^1$ with the
property that no nontrivial element of $H$ has an interval of fixed
points (e.g. $H$ is a group of real-analytic diffeomorphisms).  Then 
$H$ does not contain any finite index subgroup of:
\begin{enumerate}
\item $\Mod(g,k)$ for $g\geq 3, k\geq 0$.
\item $\Aut(F_n)$ or $\Out(F_n)$ for $n\geq 6$.
\item The Torelli group $T_{g,k}$ for $g\geq 3, k\geq 0$.
\end{enumerate}
\end{theorem}

In \cite{FF3} we construct a $C^1$ action of $T_{g,k}$ on $I$ and on
$S^1$.  Note also that $\Out(F_2)$ has a free subgroup of finite index,
which admits a faithful, $C^\omega$ action on $I$ and on $S^1$.

\bigskip
\noindent
{\bf Baumslag-Solitar groups. }The {\em Baumslag-Solitar groups}
$\BS(m,n)$ are defined by the presentation
$$\BS(m,n)=<a,b:ab^ma^{-1}=b^n>$$
When $n>m>1$ the group $\BS(m,n)$ is not residually finite; in
particular it is not a subgroup of any linear group (see, e.g.\ \cite{LS}).  We
will give a construction which shows that $\BS(m,n)$ is a subgroup of
one of the ``smallest'' infinite-dimensional Lie groups.

\begin{theorem}[Baumslag-Solitar groups: existence]
\label{theorem:bs:existence}
The group $\Diff_+^\omega(\R)$ of real-analytic diffeomorphisms of $\R$
contains a subgroup isomorphic to $\BS(m,n)$ for any $n>m\geq 1$.
The analogous result holds for $\Homeo_+(S^1)$ and $\Homeo_+(I)$.
\end{theorem}

It is not difficult to construct pairs of diffeomorphisms $a,b\in
\Diff_+^\omega(\R)$ which satisfy the relation $ab^ma^{-1}=b^n$; the
difficulty is to prove that (in certain situations) 
this is the only relation.  To do this we use a Schottky type argument.

While Ghys-Sergiescu \cite{GS} showed that $\Diff^\infty(S^1)$ contains 
Thompson's infinite simple (hence non-residually finite) group $T$, they 
also showed that $T$ admits no real-analytic action on $S^1$ (see also
\cite{FS}); indeed we do not know of any subgroups of
$\Diff^\omega(S^1)$ which are not residually finite.

The construction in the proof of Theorem \ref{theorem:bs:existence} 
gives an abundance analytic actions of $\BS(m,n)$ on $\R$, and 
$C^0$ actions of $\BS(m,n)$ on $S^1$ and on $I$.  The loss of regularity 
in moving from $\R$ to $S^1$ is no accident; 
we will show in contrast to Theorem \ref{theorem:bs:existence} 
that (for typical $m,n$) there are no $C^2$ actions of $\BS(m,n)$ 
on $S^1$ or $I$.

\begin{theorem}[Baumslag-Solitar groups: non-existence]
\label{theorem:bs:noact}
No subgroup of $\Diff_+^2(I)$ is isomorphic to $\BS(m,n),
n>m>1$.  If further $m$ does not divide $n$, then the same holds for 
$\Diff_+^2(S^1)$.
\end{theorem}

The hypothesis $m>1$ in Theorem \ref{theorem:bs:noact} is also necessary since 
$\BS(1,n),n\geq 1$ is a subgroup of $\PSL(2,\R)$, hence of
$\Diff_+^\omega(S^1)$.  

In fact $\BS(1,n)$ has many actions on $S^1$, and it is natural to
attempt a classification (up to conjugacy) of all of them.  As a first
step in this direction we observe (Theorem \ref{theorem:local_rigid}),
as a corollary of result of M. Shub \cite{Sh} on expanding maps, that
the standard, projective action of $\BS(1,n)$ on $S^1$ is {\em locally
rigid}, i.e.\ nearby actions are conjugate.  An example due to
M. Hirsch \cite{H} shows that it is not {\em globally rigid}, i.e.\
there are actions which are not conjugate to the standard action
(Theorem \ref{theorem:no_global_rigid} below).  It would be interesting to
find a numerical invariant which characterizes the standard action, as
Ghys \cite{Gh2} has done for Fuchsian groups.

\section{Tools}

In this section we recall some properties of 
diffeomorphisms of one-manifolds which will be used throughout the
paper.

\subsection{Kopell's Lemma and H\"{o}lder's Theorem}

Our primary tool is the following remarkable result of Nancy Kopell which
is Lemma 1 of \cite{K}.

\begin{thm}[Kopell's Lemma] 
\label{thm:kopell}
Suppose $f$ and $g$ are $C^2$, orientation-preserving
diffeomorphisms of an interval $[a,b)$ such that $fg = gf.$ If $f$ has no
fixed point in $(a,b)$ and $g$ has a fixed point in $(a,b)$ then
$g = id.$
\end{thm}

Another useful result is the following theorem, which is classical.

\begin{thm} [H\"older's Theorem] 
\label{thm:holder}
Suppose a group $G$ of homeomorphisms of $\R$ acts
freely and effectively on a closed subset of $\R$.  
Then $G$ is abelian.
\end{thm}

\subsection{The translation number and mean translation number}

If $f \in \Homeo_+(\T^1)$, i.e, $f$ is an orientation preserving
homeomorphism of $\T^1$, then there is a countable collection of lifts
of $f$ to orientation preserving homeomorphisms of the line.  If $F$ is
one such lift then it satisfies $FT = TF$ where $T(x) = x + 1,$ and all
others are of the form $FT^n,\ n\in \Z.$ Any orientation preserving
homeomorphism of $\R$ which commutes with $T$ is a lift of an element of
$\Homeo_+(\T^1).$ We will denote by $\Homeo_{\Z}( \R)$ the group of
homeomorphisms of $\R$ which commute with $T$, or equivalently the group
of all lifts of elements of $\Homeo_+(\T^1).$

There are two important and closely related functions from $\Homeo_+(\T^1)$ and $\Homeo_{\Z}( \R)$ to $\T^1$ 
and $\R$ respectively, which we now define.

\begin{defn}
If $F \in \Homeo_{\Z}( \R)$ define its {\em translation number} 
$\tau(F) \in \R$ by
\[
\tau(F) = \lim_{n\to \infty} \frac{F^n(x) - x}{n}
\]
If $f \in \Homeo_+(\T^1)$ define its {\em rotation number}, $\rho(f) \in \T^1,$ by
\[
\rho(f) = \tau(F) \mod (1),
\]
where $F$ is any lift of $f$.
\end{defn}

We summarize some basic properties of the rotation and translation
numbers.  Proofs of these and additional properties can be found, for
example, in \cite{dMvS}

\begin{prop}[Properties of rotation number]
\label{prop:rotation.number}
If $F \in \Homeo_{\Z}( \R)$ the number $\tau(F)$ always exists and is
independent of the $x \in \R$ used to define it.  It satisfies 
\begin{align*}
\tau(F^n) &= n \tau(F)\\
\tau(FT^n) &= \tau(F) + n\\
\tau(F_0 F F_0^{-1}) &= \tau(F) \text{ for any } F_0 \in \Homeo_{\Z}( \R).
\end{align*}
$F$ has a fixed point if and only if $\tau(F) = 0.$ 

If $f \in
\Homeo_+( \T^1)$ then $\rho(f)$ is independent of the lift $F$ used to
define it.  It satisfies
\begin{align*}
\rho(f^n) &= n \rho(f)\\
\rho(f_0 f f_0^{-1}) &= \rho(f) \text{ for any } f_0 \in \Homeo_+(\T^1).
\end{align*}
$f$ has a fixed point if and only if $\rho(f) = 0.$
\end{prop}

Suppose $G$ is a subgroup of $\Homeo_+(\T^1)$ and $\bar G$ is the group of all
lifts of elements of $G$ to $\R$.  If $G$ preserves a Borel probability 
measure $\mu_0$ on $\T^1$ then this measure may be lifted to a $\bar G$
invariant measure $\mu$ on $\R$ which is finite on compact sets and which is
preserved by the covering translation $T(x) = x +1.$
This permits us to define the {\em mean translation number} with respect to
$\mu.$

\begin{defn}
If $F \in \Homeo_{\Z}( \R)$  define its {\em mean translation number} $\tau_{\mu}(F) \in \R$ by
\begin{equation*}
\tau_{\mu}(F) = 
\begin{cases}
\mu([x, F(x)))	&\text{ if } F(x) > x,\\
0 	&\text{ if } F(x) = x,\\
-\mu([F(x), x))	&\text{ if } F(x) < x,
\end{cases}
\end{equation*}
where $x \in \R.$
\end{defn}

We enumerate some of the well known basic properties of the mean
translation number which we will use later.

\begin{prop}
\label{prop:mean.rotation}
Suppose $G$ is a subgroup of $\Homeo_+( \T^1)$ which preserves the Borel
probability measure and $\bar G$ is the group of all lifts of elements
of $G.$ If $F \in \bar G$ the mean translation number $\tau_{\mu}(F)$
is independent of the point $x$ used to define it.  Indeed
$\tau_{\mu}(F) = \tau(F).$ Moreover the function $\tau_{\mu} : \bar G
\to \R$ is a homomorphism (and hence so is $\tau : \bar G \to \R.$)
\end{prop}
\begin{proof}
Consider the function
\begin{equation*}
\nu(x,y) = 
\begin{cases}
\mu([x, y)))	&\text{ if } y > x,\\
0 	&\text{ if } y = x,\\
-\mu([y, x))	&\text{ if } y < x.
\end{cases}
\end{equation*}
It has the property that $\nu(x,y) + \nu(y,z) = \nu(x,z).$

We also note that $\tau_{\mu}(F) = \nu(x, F(x)$.  To see this is
independent of $x$ note for any $y \in \R, \ \nu(x, F(x)) = \nu(x, y)
+ \nu( y, F(y)) + \nu( F(y), F(x)) = \nu( y, F(y)) + \nu(x, y) - \nu(
F(x), F(y)).$ But $F$ preserves the measure $\mu$ and the orientation
of $\R$ so $\nu(x,y) = \nu( F(x), F(y)).$ Hence $\nu(x, F(x) = \nu( y,
F(y)).$

Since $\mu$ is the lift of a probability measure on $\T^1$
we know that $\mu([x, x+k)) = k$ for any $k \in \Z.$
So if $F^n(x) \in [x+k, x+k+1)$ we see
that 
\[
F^n(x) - x - 1 \le k \le \nu(x, F^n(x)) \le k+1 \le F^n(x) - x + 1.
\]
To see that $\tau_{\mu}(F) = \tau(F)$ we note
\begin{align*}
\tau(F) &= \lim_{n\to \infty} \frac{F^n(x) - x}{n}\\
&= \lim_{n\to \infty} \frac{\nu(x, F^n(x))}{n}\\
&=\lim_{n\to \infty}\frac{1}{n}\sum_{i=0}^{n-1}\nu(F^i(x), F^{i+1}(x))\\
&=\nu(x, F(x))\\
&=\tau_{\mu}(F).
\end{align*}
\end{proof}

There are two well known corollaries of this result which we record
here for later use.

\begin{cor}
\label{cor:rotation.homomorphism}
Suppose $G$ is a subgroup of $\Homeo_+( \T^1)$ which preserves the Borel
probability measure $\mu_0$  and $\bar G$ 
is the group of all lifts of elements of
$G.$  Then each of the functions
$\tau_{\mu} : \bar G \to \R, \ \tau : \bar G \to \R$ and 
$\rho : G \to \T^1$ is a homomorphism.
\end{cor}

\begin{proof}
To see that $\tau_{\mu} : \bar G \to \R$ is a homomorphism we suppose
$f,g \in \bar G$ and consider $\nu$ as defined above.  Then
$\tau_{\mu}(fg) = \nu( x, f(g(x))) = \nu( x, g(x)) + \nu(
g(x),f(g(x))) = = \tau_{\mu}(g) + \tau_{\mu}(f).$

Since $\tau_{\mu} = \tau$ we know that $\tau$ is also a homomorphism.
There is a natural homomorphism $\pi: \bar G \to G$ which assigns to
$g$ its projection on $\T^1$.  
If $p: \R \to \T^1 = \R/\Z$ is the natural projection then
for any $f,g \in \bar G,$ 
\[
\rho( \pi(fg)) =  p(\tau(fg)) = p(\tau(f)) + p(\tau(g)) = \rho( \pi(f)) + \rho( \pi(g)).
\]
Hence $\rho: G \to \T^1$ is a homomorphism.
\end{proof}

\begin{cor}
\label{cor:abelian.subgroup}
If $G$ is an abelian subgroup of $\Homeo_+(\T^1)$ and $\bar G$ is the
group of all lifts to $\R$ of elements of $G$, then $\bar G$
is abelian and both $\tau: \bar G \to \R$ and 
$\rho: G \to \T^1$ are homomorphisms.
\end{cor}
\begin{proof}
Since $G$ is abelian it is amenable and there is a Borel probability
measure $\mu_0$  on $\T^1$ invariant under $G$.  Let $\mu$ be the
lift of this measure to $\R.$  Then Corollary \ref{cor:rotation.homomorphism}
implies that $\tau : \bar G \to \R$ and $\rho: G \to \T^1$ are homomorphisms.

As above let $\pi: \bar G \to G$ be the natural projection.  The
kernel of $\pi$ is $\{T^n\}_{n \in \Z}$, where $T(x) = x+1$, i.e., the
lifts of the identity. If $f,g \in \bar G$ then $[f,g]$ is in the
kernel of $\pi$ so $[f,g] = T^k,$ for some $k \in \Z.$ But since $\bar
\rho$ is a homomorphism $\rho([f,g]) = 0$ which implies $k=0.$ Hence
$\bar G$ is abelian.
\end{proof}

The following easy lemma turns out to be very useful.

\begin{lem}
\label{lem:equal.on.support}
Suppose $G$ is a subgroup of $\Homeo_+( \T^1)$ which preserves the Borel
probability measure $\mu_0$  and $f_0,g_0 \in G$  satisfy 
$\rho(f_0) = \rho(g_0).$   Then $f_0(x) =  g_0(x)$ for
all $x$ in the support of $\mu_0.$
\end{lem}
\begin{proof}
Let $\mu$ be the lift of the measure $\mu_0$ to $\R.$ Pick lifts $f$
and $g$ of $f_0$ and $g_0$ respectively which satisfy $\tau_{\mu}( f)
= \tau_{\mu}( g).$ Then for any $x \in \R$ we have $\mu([x,f(x))) =
\tau_{\mu}( f) = \tau_{\mu}( g) = \mu([x,g(x))).$

Suppose that $g(x) > f(x)$.  It then follows that $\mu([f(x),g(x)) =
0,$ so for a sufficiently small $\epsilon >0, \
\mu([f(x),f(x)+\epsilon])= 0$ and $\mu([g(x)-\epsilon, g(x)))= 0.$
Applying $gf^{-1}$ to $[f(x),f(x)+\epsilon)$ we see that $\mu([g(x),
g(x)+\epsilon'))= 0$.  Hence $\mu([g(x)-\epsilon, g(x)+\epsilon'))= 0$
which is not possible if $x$ and hence $g(x)$ is in the support of
$\mu.$

We have shown that $x \in \supp(\mu)$ implies $f(x) \le g(x).$
The inequality $g(x) \le f(x)$ is similar.  Projecting back to
$\T^1$ we conclude that $f_0 = g_0$ on $\supp(\mu_0).$
\end{proof}

One additional well known property we will use is the following.

\begin{prop}
\label{prop:irrational.rotation}
If $g \in \Homeo_+(\T^1)$ is an irrational rotation then its centralizer 
$Z(g)$ in $\Homeo_+(\T^1)$ is the group of rigid rotations of $\T^1.$
\end{prop}
\begin{proof}
Let $f$ be an element of $Z(g)$.  Then the group generated by $f$ and $g$
is abelian and hence amenable so there it has an invariant Borel probability
measure.  But Lebesgue measure is the unique Borel probability measure
invariant by the irrational rotation $g$.  Since $f$ preserves orientation
and  Lebesgue measure, it is a rotation.
\end{proof}

\section{Mapping Class Groups}

In this section we prove Theorem \ref{theorem:mcgroups}. Several of the
results, in particular those in 
\S\ref{subsection:fully:supported}, will also be used to prove 
Theorem \ref{theorem:autogroups}.
	
\subsection{Fully Supported Diffeomorphisms}
\label{subsection:fully:supported}

One of the main techniques of this paper is to use Kopell's Lemma
(Theorem \ref{thm:kopell}) to understand actions of commuting
diffeomorphisms.  The key dichotomy that arises is the behavior of 
diffeomorphisms with an interval of fixed points versus those which we
call {\em fully supported}.

If $f$ is a homeomorphism of a manifold $M$, we will denote by
$\partial\Fix(f)$ the frontier of $\Fix(f)$, i.e. the set
$\partial\Fix(f) =\Fix(f) \setminus \Int(\Fix(f)).$

\begin{defn}[Fully supported homeomorphism] 
A homeomorphism $f$ of a manifold $M$ is
{\em fully supported} provided that $\Int(\Fix(f)) = \emptyset$, or
equivalently $\partial\Fix(f) =\Fix(f).$  
A subgroup $G$ of $\Homeo(M)$
will be called {\em fully supported} provided that every nontrivial
element is fully supported.
\end{defn} 

It is a trivial but useful observation that if
$M$ is connected then for every nontrivial homeomorphism $f$ of $M$, the
set $\Fix(f) \ne \emptyset$ if and only if $\partial \Fix(f) \ne
\emptyset.$

The following lemma indicates a strong consequence of Kopell's Lemma (Theorem 
\ref{thm:kopell}).

\begin{lem}[Commuting diffeomorphisms on $I$]
\label{lem:fix.component.I}
Suppose $f$ and $g$ are commuting orientation preserving $C^2$
diffeomorphisms of $I$.  Then $f$ preserves every component of
$\Fix(g)$ and vice versa.  Moreover, $\partial \Fix(f)) \subset
\Fix(g)$ and and vice versa.  In particular if $f$ and $g$ are fully
supported then $\Fix(f) = \Fix(g)$.
\end{lem}

\begin{proof}
The proof is by contradiction.  Assume $X$ is a component of $\Fix(g)$
and $f(X) \ne X.$ Since $f$ and $g$ commute $f(\Fix(g)) = \Fix(g)$ so
$f(X) \ne X$ implies $f(X) \cap X = \emptyset.$ Let $x$ be an element
of $X$ and without loss of generality assume $f(x) < x.$  Define
\[
a = \lim_{n \to \infty} f^n(x) \text{ and } b = \lim_{n \to -\infty} f^n(x).
\]
Then $a$ and $b$ are fixed under both $f$ and $g$ and $f$ has no fixed
points in $(a,b).$ Then Kopell's Lemma (Theorem \ref{thm:kopell})
implies $g(y) = y$ for all $y \in [a,b]$, contradicting the hypothesis
that $X$ is a component of $\Fix(g)$.  The observation that $\partial
\Fix(f) \subset \Fix(g)$ follows from the fact that $x \in \partial
\Fix(f)$ implies that either $\{x\}$ is a component of $\Fix(f)$ or
$x$ is the endpoint of an interval which is a component of $\Fix(f)$ so,
in either case, $x \in \Fix(g).$
\end{proof}

There is also a version of Lemma \ref{lem:fix.component.I} for the
circle.  For $g\in \Homeo_+(\T^1)$ let $\Per(g)$ denote the set of
periodic points of $g$.

\begin{lem}[Commuting diffeomorphisms on $\T^1$]
\label{lem:per.component}
Suppose $f$ and $g$ are commuting orientation preserving $C^2$
diffeomorphisms of $\T^1.$ Then $f$ preserves every component of
$\Per(g)$ and vice versa.  Moreover, $\partial \Per(f) \subset
\Per(g)$ and vice versa.  In particular if neither $\Per(f)$ or
$\Per(g)$ has interior then $\Per(f) = \Per(g)$.
\end{lem}

\begin{proof}
First consider the case that $\Per(g) = \emptyset.$ The assertion that
$f$ preserves components of $\Per(g)$ is then trivial.  Also, $g$ is
topologically conjugate to an irrational rotation which by Proposition
\ref{prop:irrational.rotation} implies that $f$ is topologically
conjugate to a rotation and hence $\Per(f) = \emptyset$ or $\Per(f) =
\T^1.$ In either case we have the desired result.

Thus we may assume both $f$ and $g$ have periodic points.  Since for
any circle homeomorphism all periodic points must have the same period, 
we can let $p$ be the least common multiple of the periods and observe
that $\Per(f) = \Fix(f^p)$ and $\Per(g) = \Fix(g^p).$ Since $f$ and
$g$ commute and both $f^p$ and $g^p$ have fixed points, if $x \in
\Fix(f^p)$ then $y =\lim_{n \to \infty}g^{np}(x)$ will exist and be a
common fixed point for $f^p$ and $g^p$.  If we split the circle at $y$
we obtain two commuting diffeomorphisms, $f^p$ and $g^p,$ of an
interval to which we may apply Lemma \ref{lem:fix.component.I} and
obtain the desired result.
\end{proof}

\begin{lem}
\label{lem:rho.centralizer}
Let $g_0 \in \Diff_+^2(\T^1)$ and let $Z(g_0)$ denote its centralizer.
Then the rotation number $\rho : Z(g_0) \to \T^1$ is a homomorphism.
\end{lem}
\begin{proof}
We first observe the result is easy if $g_0$ has no periodic points.
In this case $g_0$ is conjugate in $\Homeo_+(\T^1)$ to an irrational
rotation by Denjoy's Theorem, and the centralizer of an irrational
rotation is abelian by Proposition \ref{prop:irrational.rotation}.  So
the centralizer of $g$ in $\Diff_+^1(\T^1)$ is abelian.  But then $\rho$
restricted to an abelian subgroup of $\Homeo_+(\T^1)$ is a homomorphism
by Proposition \ref{cor:abelian.subgroup}.

Thus we may assume $g_0$ has a periodic point, say of period $p$.
Then $h_0 = g_0^p$ has a fixed point. Let $h$ be a lift of $h_0$ to
$\R$ which has a fixed point.  Let $G$ denote the group of all lifts
to $\R$ of all elements of $Z(g_0)$.  Note that every element of this
group commutes with $h$ because by Corollary
\ref{cor:abelian.subgroup} any lifts to $\R$ of commuting
homeomorphisms of $\T^1$ commute.

Let $X = \partial \Fix(h)$ and observe that $g(X) = X$ for all $g\in
G$, so $G$ acts on the unbounded closed set $X$ and $G$ acts on $X_0 =
\partial \Fix(h_0)$ which is the image of $X$ under the covering
projection to $\T^1$ .  We will show that if $\Fix(g) \cap X \ne
\emptyset$ then $X \subset \Fix(g).$ Indeed applying Lemma
\ref{lem:per.component} to $h_0$ and map $g_1$ of $\T^1$ which $g$
covers, we observe that $\partial \Fix(h) \subset \Fix(g_1)$ which
implies $X \subset \partial \Fix(g)$ if $g$ has has a fixed point.

It follows that if $H = \{ g \in G \ | g(x) = x \text{ for all } x \in
X\}$ is the stabilizer of $X$ then $G/H$ acts freely on $X$ and hence
is abelian by H\"older's theorem.  Also $G/H$ acts on $X_0$ and hence
there is a measure $\mu_0$ supported on $X_0$ and invariant under $G/H$.
Clearly this measure is also invariant under the action of $Z(g_0)$ on
$\T^1$.  The measure $\mu_0$ lifts to a $G$-invariant measure $\mu.$

It follows from Corollary \ref{cor:rotation.homomorphism}
that the translation number $\tau: G \to \R$ 
and the rotation number $\rho: Z(g_0) \to \T^1$  are homomorphisms.
\end{proof}

The {\em commutativity graph} of a set $S$ of generators of a group is the
graph consisting of one vertex for each $s\in S$ and an edge connecting 
elements of $S$ which commute.  

\begin{thm}[abelian criterion for $I$]
\label{thm:comm.graph.I}
Let $\{g_1,g_2, \dots, g_k\}$ be a set of fully supported elements of
$\Diff_+^2(I)$ and let $G$ be the group they generate.  Suppose that the 
commutativity graph of this generating set is connected. Then $G$ is
abelian.
\end{thm}
\begin{proof}
By Lemma \ref{lem:fix.component.I} we may conclude that for each
$j_1,j_2$ we have $\Fix(g_{j_1}) = \Fix(g_{j_2}).$ Call
this set of fixed points $F$.  Clearly $F$ is the set of global fixed
points of $G$.

Fix a value of $i$ and consider $Z(g_i).$ Restricting to any component
$U$ of the complement of $F$ we consider the possiblity that there is
an $h \in Z(g_i)$ with a fixed point in $U$.  Kopell's Lemma 
(Theorem \ref{thm:kopell}), applied to the closure of the open inteval $U$,
tells us that such an $h$ is the identity on $U$.
Thus the restriction of $Z(g_i)$ to $U$ is free and hence abelian by
H\"older's Theorem.  But $U$ was an arbitrary component of the
complement of $F$ and obviously elements of $Z(g_i)$ commute on their
common fixed set $F$.  So we conclude that $Z(g_i)$ is abelian.

We have hence shown that if $g_j$ and $g_k$ are joined by a path of
length two in the commutativity graph, they are joined by a path of 
length one.  A straightforward induction shows that any two generators
are joined by a path of length one, i.e. any two commute.
\end{proof}

There is also a version of Theorem \ref{thm:comm.graph.I} for the
circle.

\begin{cor}[Abelian criterion for $\T^1$]
\label{cor:comm.graph.T^1}
Let $\{g_1,g_2, \dots, g_k\}$ be a set of fully supported elements of
$\Diff_+^2(\T^1)$, each of which has a fixed point, and let $G$ be the
group they generate.  Suppose that the commutativity graph of
these generators is connected. Then $G$ is abelian.
\end{cor}

\begin{proof}
By Lemma \ref{lem:per.component} we may conclude that for each $j,k$
we have $\Per(g_j) = \Per(g_k).$ But since each $g_i$ has a fixed
point $\Per(g_i) = \Fix(g_i).$  Hence there is a common fixed point for
all of the generators. Splitting $\T^1$ at a common fixed point we see
$G$ is isomorphic to a subgroup of $\Diff_+^2(I)$ satisfying the
hypothesis of Theorem \ref{thm:comm.graph.I}.  It follows that
$G$ is abelian.
\end{proof}

\noindent
{\bf Question. }  Let $G$ be a subgroup of $\Homeo_+(\T^1$).  
If the commutativity graph for a set of generators of
$G$ is connected, does this imply that rotation
number is a homomorphism when restricted to $G$?
\bigskip

We will need to understand relations between fixed sets not just
of commuting diffeomorphisms, but also of diffeomorphisms with
another basic relation which occurs in mapping class groups as
well as in automorphism groups of free groups.

\begin{lem}[\boldmath$aba$ lemma]
\label{lem:aba}
Suppose $a$ and $b$ are elements of $\Homeo_+(I)$ or elements of
$\Homeo_+(\T^1)$ which have fixed points.  Suppose also that $a$
and $b$ satisfy the relation $a^{n_1}b^{m_3}a^{n_2} =
b^{m_1}a^{n_3}b^{m_2}$ with $m_1 + m_2 \ne m_3$.  If $J$ is a
nontrivial interval in $\Fix(a)$ then either $J \subset \Fix(b)$
or $J \cap \Fix(b) = \emptyset.$
\end{lem}
\begin{proof}
Suppose $J \subset \Fix(a)$ and suppose $z \in J \cap \Fix(b).$ We
will show that this implies that $J \subset \Fix(b).$

Let $x_0 \in J.$ Observe that $a^{n_1}b^{m_3}a^{n_2} =
b^{m_1}a^{n_3}b^{m_2}$ implies $a^{-n_2}b^{-m_3}a^{-n_1} =
b^{-m_2}a^{-n_3}b^{-m_1}$.  Hence we may assume without loss of
generality that $b^{m_3}(x_0)$ is in the subinterval of $J$ with
endpoints $x_0$ and $z$ (otherwise replace $a$ and $b$ by their
inverses).  In particular $b^{m_3}(x_0) \in J.$ Note that 
$$b^{m_3}(x_0) =a^{n_1}b^{m_3}(x_0) = a^{n_1}b^{m_3}a^{n_2}(x_0) =
b^{m_1}a^{n_3}b^{m_2}(x_0) = b^{m_1}b^{m_2}(x_0) = b^{m_1+m_2}(x_0)$$
This implies $b^{m_1+m_2 -m_3}(x_0) = x_0$ and since $m_1+m_2 -m_3 \ne
0$ we can conclude that $b(x_0) = x_0.$  Since $x_0 \in J$ was arbitrary
we see $J \subset \Fix(b).$
\end{proof}

\subsection{Groups of type \boldmath$MC(n)$}

We will need to abstract some of the properties of the
standard generating set for the mapping class group of a surface in
order to apply them in other circumstances.

\begin{defn}[Groups of type \boldmath$MC(n)$]
 \label{defn:mc(n)}
We will say that a group $G$ is {\em of type $MC(n)$} provided it is
nonabelian and has a set of generators $\{ a_i,b_i\}_{i=1}^n \cup
\{c_j\}_{j=1}^{n-1}$ with the properties that
\begin{enumerate}
\item $a_i a_j = a_j a_i,\ b_i b_j = b_j b_i,\ c_i c_j = c_j c_i,\ a_i
c_j = c_j a_i,$ for all $i,j,$
\item $a_i b_j = b_j a_i,$ if $i \ne j$,
\item $b_i c_j = c_j b_i,$ if $j \ne i, i-1$,
\item $a_ib_ia_i = b_i a_i b_i$ for $1 \le i \le n,$ and
$b_i c_j b_i = c_jb_ic_j$ whenever $j = i, i-1$.
\end{enumerate}
\end{defn}

A useful consequence of the $MC(n)$ condition is the following.

\begin{lem}
\label{lemma:conj}
If $G$ is a group of type $MC(n)$ and $\{ a_i,b_i\}_{i=1}^n \cup
\{c_j\}_{j=1}^{n-1}$ are the generators guaranteed by the definition
then any two of these generators are conjugate in $G$.
\end{lem}
\begin{proof}
The relation $aba = bab$ implies $a = (ab)^{-1}b(ab)$, so
$a$ and $b$ are conjugate.  Therefore property (4) of the definition
implies $a_i$ is conjugate to $b_i$ and $b_i$ is conjugate to 
$c_i$ and $c_{i-1}$. This proves the result.
\end{proof}

The paradigm of a group of type $MC(n)$ is the mapping class
group $\Mod(n,0)$.

\begin{prop}[\boldmath$\Mod(n,0)$ and \boldmath$\Mod(n,1)$]
\label{lem:conj}
For $n>2$ and $k=0,1$ the group $\Mod(n,k)$ contains a set of elements 
$${\cal S}=\{ a_i, b_i\}_{i=1}^n \cup \{c_j\}_{j=1}^{n-1}$$ with the following
properties:
\begin{enumerate}
\item For $k=0$ the set $\cal S$ generates all of $\Mod(n,0)$.
\item The group generated by $\cal S$ is of type $MC(n)$.
\item There is an element of $g\in \Mod(n,k)$ 
such that $g^{-1}a_1g = a_n,\ g^{-1}b_1g = b_n,\ g^{-1}b_2g = b_{n-1},$
and $g^{-1}c_1g = c_{n-1}.$
\item If $G_0$ is the subgroup of $\Mod(n,k)$ generated by the 
subset  $$\{ a_i, b_i\}_{i=1}^{n-1} \cup \{c_j\}_{j=1}^{n-2}$$ 
then $G_0$ is of type $MC(n-1)$.
\end{enumerate}
\end{prop}
\begin{proof}
Let $\Sigma$ denote the surface of genus $n$ with $k$ punctures.  
Let $i(\alpha,\beta)$ denote the geometric intersection number of the 
(isotopy classes of) closed curves $\alpha$ and $\beta$ on $\Sigma$.    
Let $a_i,b_i,c_j$ with $1\leq i\leq n, 1\leq j\leq n$ be Dehn
twists about essential, simple closed curves
$\alpha_i,\beta_i,\gamma_j$ with the following properties:
\begin{itemize}
\item
$i(\alpha_i,\alpha_j)=i(\beta_i,\beta_j)=i(\gamma_i,\gamma_j)=0$
for $i\neq j$.

\item $i(\alpha_i,\beta_i)=1$ for each $1\leq i\leq n$.
\item $i(\gamma_j,\beta_{j})=1=i(\gamma_j,\beta_{j+1})$ for
$1\leq j\leq n-1$.
\end{itemize} 

It is known (see, e.g.\cite{Iv}) that there exist such
$a_i,b_i,c_j$ which generate $\Mod(n,0)$; we choose such
elements.  These are also nontrivial in $\Mod(n,1)$.  
Since Dehn twists about simple closed curves with
intersection number zero commute, and since Dehn twists $a,b$
about essential, simple closed curves with intersection number
one satisy $aba=bab$ (see, e.g.\ Lemma 4.1.F of \cite{Iv}), 
properties (1)-(4)  in Definition \ref{defn:mc(n)} follow.  We note that 
the proofs of these relations do not depend on the location of the
puncture, as long as it is chosen off the curves $\alpha_i,\beta_i,\gamma_j$.

To prove item (2) of the lemma, consider the surface of genus $g-2$ and
with $2$ boundary components obtained by cutting $\Sigma$ along $\alpha_1$
and $\gamma_1$.  The loop $\beta_1$ becomes a pair of arcs connecting
$2$ pairs of boundary components, and $\beta_2$ becomes an arc
connecting $2$ boundary components.  The genus, boundary components, and 
combinatorics of arcs is the exact same when cutting $\Sigma$ along 
$\alpha_n$ and $\gamma_{n-1}$.  Hence by the classification of surfaces
it follows that there is a homeomorphism between resulting surfaces,
inducing a homeomorphism $h:\Sigma\rightarrow \Sigma$ with 
$h(\alpha_1)=\alpha_n, h(\gamma_1)=\gamma_{n-1}$ and 
$h(\beta_i)=\beta_{n-i+1}, i=1,2$.  
The homotopy class of $h$ gives the required 
element $g$.

To prove item (3), let $\tau$ be an essential, separating, simple closed
curve on $\Sigma$ such that one of the components $\Sigma'$ of
$\Sigma-\tau$ has genus one and contains $a_n$ and $b_n$.  Then there is
a homeomorphism $h\in\Homeo_+(\Sigma)$ taking any element of 
$\{ a_i, b_i\}_{i=1}^{n-1} \cup \{c_j\}_{j=1}^{n-2}$ to any other
element, and which is the identity on $\Sigma'$.  Then the isotopy class 
of $h$, as an element of $\Mod(n,0)$, is the required conjugate, since
it lies in the subgroup of $\Mod(0,n)$ of diffeomorphisms supported on 
$\Sigma-\Sigma'$, which is isomorphic to $\Mod(n-1,1)$ and equals $G_0$.
\end{proof}

\subsection{Actions on the interval}

In this subsetion we consider actions on the interval $I=[0,1]$.

\begin{thm}[\boldmath$MC(n)$ actions on $I$]
\label{thm:mc(n)}
Any $C^2$ action of a group $G$  of type $MC(n)$ for $n \ge 2$ on
an interval $I$ is abelian.
\end{thm}
\begin{proof}

We may choose a set of generators $\{ a_i, b_i\}_{i=1}^n \cup
\{c_j\}_{j=1}^{n-1}$ for $G$ with the properties listed in Definition
\ref{defn:mc(n)}.

If $G$ has global fixed points other than the endpoints of $I$ we wish
to consider the restriction of the action to the closure of a component
of the complement of these global fixed points.  If all of these restricted
actions are abelian then the original action was abelian.  Hence it suffices
to prove that these restrictions are abelian.
Thus we may consider an action of $G$ on a closed interval $I_0$
with no interior global fixed points.  None of the generators above
can act trivially on $I_0$ since the fact that they are all conjugate
would mean they all act trivially.  

We first consider the case that one generator has a nontrivial
interval of fixed points in $I_0$ and show this leads to a
contradiction.  Since they are all conjugate, each of the generators
has a nontrivial interval of fixed points.  

Choose an interval $J$ of fixed points which is maximal among all
intervals of fixed points for all of the $a_i$.  That is, $J$ is a
nontrivial interval of fixed points for one of the $a_i$, which we
assume (without loss of generality) is $a_2$, and there is no $J'$ which
properly contains $J$ and which is pointwise fixed by some $a_j$.

Suppose $J = [x_0, x_1]$.  At least one of $x_0$ and $x_1$ is not an
endpoint of $I_0$ since $a_2$ is not the identity.
Suppose it is $x_0$ which is in
$\Int(I_0).$ Then $x_0 \in \partial \Fix(a_2)$, so by
Lemma \ref{lem:fix.component.I} we know that $x_0$ is fixed by all of
the generators except possibly $b_2$, since $b_2$ is the only generator
with which $a_2$ does not commute.  The point $x_0$ cannot be fixed by
$b_2$ since otherwise it would be an interior point of $I_0$ fixed by all the
generators, but there are no global fixed points in $I_0$ other than the
endpoints.  The identity $a_2b_2a_2 = b_2a_2b_2$ together with Lemma
\ref{lem:aba} then tells us that $b_2(J) \cap J = \emptyset.$

Now assume without loss of generality that $b_2(x_0) > x_0$ (otherwise
replace all generators by their inverses). 
Define 
\[
y_0 = \lim_{n \to -\infty} b_2^n(x_0) \text{ and } y_1 = \lim_{n \to
\infty} b_2^n(x_0).
\]
Since $b_2(J) \cap J = \emptyset$ we have $b_2^{-1}(J) \cap J =
\emptyset$ and hence $y_0 < x_0 < x_1 < y_1.$ But $x_0$ is a fixed
point of $a_1$ (as well as the other generators which commute with
$a_2$) and $a_1$ commutes with $b_2.$ Hence $y_0$ and $y_1$ are fixed
by $a_1$ (since $b_2$ preserves $\Fix(a_1)$). We can now apply
Kopell's Lemma (Theorem \ref{thm:kopell}) 
to conclude that $a_1$ is the identity
on $[y_0, y_1]$ which contradicts the fact that $J = [x_0,x_1]$ is a
maximal interval of fixed points among all the $a_i$.

We have thus contradicted the supposition that one of the generators
has a nontrivial interval of fixed points in $I_0$, so we may assume
that each of the generators, when restricted to $I_0$, has fixed point
set with empty interior.  That is, each generator is fully supported
on $I_0.$ 

We now note that given any two of the generators above, there is
another generator $h$ with which they commute.  Thus we may conclude
from Lemma \ref{thm:comm.graph.I} that the action of $G$ on $I_0$
abelian.
\end{proof}

\noindent
{\bf Proof of Theorem \ref{theorem:mcgroups} for $I$: }
By Proposition \ref{lem:conj}, the mapping class group $\Mod(g,k), g\geq
2, k=0,1$ is a group of type $MC(g).$ Since every abelian quotient of
$\Mod(g,k), g\geq 2$ is finite (see, e.g.\ \cite{Iv}), in fact trivial
for $g\geq 3$, and since finite groups must act trivially on $I$, the
statement of Theorem \ref{theorem:mcgroups} for $I$ follows.

\subsection{Actions on the circle}

We are now prepared to prove Theorem \ref{theorem:mcgroups} in the case
of the circle $\T^1$.

\bigskip
\noindent
{\bf Theorem \ref{theorem:mcgroups} (Circle case). }
{\it Any $C^2$ action of the mapping class group 
$M(g,k)$  of a surface of genus $g \ge 3$
with $k = 0, 1$ punctures on $\T^1$ must be trivial.}
\bigskip

As mentioned in the introduction, the $C^2$ hypothesis is necessary.

\begin{proof}

\medskip
\noindent
{\bf Closed case. }
We first consider the group $G=\Mod(g,0)$, and choose 
generators $\{ a_i, b_i\}_{i=1}^n \cup
\{c_j\}_{j=1}^{n-1}$ for $G$ with the properties listed in Definition
\ref{defn:mc(n)} and Proposition \ref{lem:conj}.  
All of these elements are conjugate by Lemma \ref{lemma:conj}, so they all
have the same rotation number.

We wish to consider a subgroup $G_0$ of type $MC(n-1).$  We let
$u = a_n^{-1}$ and for $1 \le i \le n-1$ we set $A_i = a_i u$ and 
$\ B_i = b_i u$.  For $1 \le i \le n-2$ we let $C_i = c_i u.$
Since $u$ commutes with any element of $\{ a_i,b_i\}_{i=1}^{n-1} \cup \{c_j\}_{j=1}^{n-2}$
we have the relations
\begin{enumerate}
\item $A_i A_j = A_j A_i,\ B_i B_j = B_j B_i,\ C_i C_j = C_j C_i,\ A_i
C_j = C_j A_i,$ for all $i,j,$
\item $A_i B_j = B_j A_i,$ if $i \ne j$,
\item $B_i C_j = C_j B_i,$ if $j \ne i, i+1$,
\item $A_iB_iA_i = B_iA_iB_i$ and
$B_i C_j B_i = C_j B_i C_j$ whenever $j = i, i-1$.
\end{enumerate}

We now define $G_0$ to be the subgroup of $G$ generated by 
$\{A_i,B_i\}_{i=1}^{n-1} \cup \{C_j\}_{j=1}^{n-2}$.  The fact that $u$
commutes with each of $\{ a_i,b_i\}_{i=1}^{n-1} \cup
\{c_j\}_{j=1}^{n-2}$ and has the opposite rotation number implies that
the rotation number of every element of $\{ A_i,B_i\}_{i=1}^{n-1}
\cup \{C_j\}_{j=1}^{n-2}$ is $0.$ Thus each of these elements has a
fixed point.

We note that given any two of these generators of $G_0$ there is
another generator with which they commute.  If all of these generators
have fixed point sets with empty interior then we may conclude from
Lemma \ref{lem:fix.component.I} that any two of them have equal fixed
point sets, i.e. that $\Fix(A_i) = \Fix(B_j) = \Fix(C_k)$ for all
$i,j,k.$ So in this case we have found a common fixed point for all
generators of $G_0.$ If we split $\T^1$ at this point we get an action
of $G_0$ on an interval $I$ which must be abelian by Theorem
\ref{thm:mc(n)}.  It follows that all of the original generators
$\{ a_i, b_i\}_{i=1}^n \cup \{c_j\}_{j=1}^{n-1}$ commute with each other 
except possibly $c_{n-1}$ may not commute with $b_{n-1}$ and
$b_n$, and $a_n$ and $b_n$ may not commute.  

Let $\psi:G\rightarrow \Diff_+^2(\T^1)$ be the putative action, 
and consider the element
$g$ guaranteed by Proposition \ref{lem:conj}.  We have that 
$$
[\psi(a_n),\psi(b_n)] =\psi(g)^{-1}[\psi(a_1),\psi(b_1)]\psi(g) 
=id$$
and
$$
[\psi(c_{n-1}),\psi(b_n)] = 
\psi(g)^{-1}[\psi(c_1),\psi(b_1)]\psi(g)
=id 
$$
and
$$[\psi(c_{n-1}),\psi(b_{n-1})] =\psi(g)^{-1}[\psi(c_1),\psi(b_2)]\psi(g)
=id
$$
Hence $G$ is abelian, hence trivial since (see \cite{Iv}) abelian 
quotients of $\Mod(g,k), g\geq 3$ are trivial.

Thus we are left with the case that one generator of $G_0$ has a
nontrivial interval of fixed points in $\T^1.$ Since they are all
conjugate, each of the generators of $G_0$ 
has has a nontrivial interval of fixed
points.  Also we may assume no generator fixes every point of
$\T^1$ since if one did the fact that they are all conjugate would imply they all
act trivially.

Choose a maximal interval of fixed points $J$ for any of the subset of
generators $\{A_i\}.$ That is, $J$ is a nontrivial interval of fixed
points for one of the $\{A_i\}$, which we assume (without loss of
generality) is $A_2$, and there is no $J'$ which properly contains $J$
and which is pointwise fixed by some $A_j$.

Suppose $x_0$ and $x_1$ are the endpoints of $J$.  Then $x_0,x_1 \in 
\partial \Fix(A_2)$
so by Lemma \ref{lem:per.component}, $x_0$ and $x_1$
are fixed points for all of the generators except $B_2$, since
$B_2$ is the only one with which $A_2$ does not commute.  If the point
$x_0 \in \Fix(B_2)$ we have found a common fixed point for all the
generators of $G_0$ and we can split $\T^1$ at this point obtaining an
action of $G_0$ on an interval which implies $G_0$ is abelian by
Theorem \ref{thm:mc(n)}.  So suppose $x_0 \notin \Fix(B_2)$.

Recall that the diffeomorphism $B_2$ has rotation number $0$, so it
fixes some point.  Thus the identity $A_2B_2A_2 = B_2A_2B_2$, together
with the fact that $B_2$ cannot fix $x_0$, implies by Lemma
\ref{lem:aba} that $B_2(J) \cap J = \emptyset$ and $J \cap B_2^{-1}(J) =
\emptyset.$ Define
\[
y_0 = \lim_{n \to -\infty} B_2^n(x_0) \text{ and } y_1 = \lim_{n \to
\infty} B_2^n(x_0).
\]
We will denote by $K$ the interval in $\T^1$ with endpoints $y_0$ and $y_1$
which contains $x_0$.  Note that $K$ properly contains $J$ and that 
$B_2(K) = K.$

Since $x_0 \in \Fix(A_1)$ and $A_1$ commutes with $B_2$ we conclude
$y_0$ and $y_1$ are fixed by $A_1$. We can now apply Kopell's Lemma 
(Theorem \ref{thm:kopell}) to the interval $K$ with diffeomorphisms $B_2$ and
$A_1$ to conclude that $A_1$ is the identity on $K$.  But $J$ is a
proper subinterval of $K$ which contradicts the fact that $J$ was a
maximal interval of fixed points for any one of the $A_i$.  Thus in this
case too we have arrived at a contradiction.

\bigskip
\noindent
{\bf Punctured case. }We now consider the group $\Mod(g,1)$, and
choose elements $\{ a_i, b_i\}_{i=1}^n \cup
\{c_j\}_{j=1}^{n-1}$ as above.  The argument above shows that the
subgroup $G$ of $\Mod(g,1)$ generated by these elements 
\footnote{We do not know wether or not $G$ actually equals $\Mod(g,1)$;
there seems to be some confusion about this point in the literature.}
acts trivially on $S^1$.  We claim that the normal closure of $G$ in
$\Mod(g,1)$ is all of $\Mod(g,1)$, from which it follows that
$\Mod(g,1)$ acts trivially, finishing the proof.

To prove the claim, let $S_g$ denote the closed surface of genus $g$ 
and recall the exact sequence (see, e.g.\ \cite{Iv}):
$$1\rightarrow \pi_1(S_g)\rightarrow \Mod(g,1)\rightarrow
\Mod(g,0)\rightarrow 1$$
where $\pi_1(Sigma_g)$ is generated by finitely many ``pushing the
point'' homeomorphisms $p_\tau$ around generating loops $\tau$ 
in $\pi_1(S_g)$ with
basepoint the puncture.  Each generating loop $\tau$ has intersection number
one with exactly one of the loops, say $\beta_i$, corresponding to one
of the twist generators of $G$.  Conjugating $b_i$ by $p_\tau$ 
gives a twist about a loop $\beta'_i$ which together with
$\beta_i$ bounds an annulus containing the puncture.  The twist about
$\beta_i$ composed with a negative twist about $\beta'_i$ gives the
isotopy class of $p_\tau$.  In this way we see that the normal closure
of $G$ in $\Mod(g,1)$ contains $G$ together with each of the generators
of the kernel of the above exact sequence, proving the claim.

\end{proof}

\section{$\Aut(F_n)$, $\Out(F_n)$ and other subgroups}

In \S\ref{section:aut1} we prove Theorem \ref{theorem:autogroups}.  Note
that the result of Bridson-Vogtmann mentioned in the introduction
implies this theorem since it is easy to see that any finite subgroup of
$\Homeo_+(\T^1)$ is abelian so their result implies ours.  We give our
proof here because it is short, straightforward and provides a good
illustration of the use of the techniques developed above.  While
several aspects of the proof are similar to the case of $\Mod(g,0)$, we
have not been able to find a single theorem from which both results
follow.

In \S\ref{section:aut2} we prove theorem \ref{theorem:finite:index},
extending the application to finite index subgroups of $\Mod(g,k)$ and
$\Aut(F_n)$.  

\subsection{Actions of $\Aut(F_n)$ and $\Out(F_n)$}
\label{section:aut1}
We begin with a statement of a few standard facts about generators and
relations of $\Aut( F_n)$.

\begin{lem}[Generators for \boldmath$\Aut(F_n)$]
\label{lem:generators}
The group $\Aut( F_n)$ has a subgroup of index two which has a 
set of generators $\{A_{ij}, B_{ij}\}$ with
$i \ne j,\ 1 \le i \le n$ and $1 \le j \le n.$ These generators satisfy
the relations
\begin{gather}
A_{ij} A_{kl} = A_{kl} A_{ij} \text{ and } B_{ij} B_{kl} = B_{kl} B_{ij}
\text{ if }
\{i,j\} \cap \{k,l\} = \emptyset \tag{i}\\
[A_{ij}, A_{jk}] = A^{-1}_{ik} \text{ and }[A_{ij}, A^{-1}_{jk}] =
A_{ik}\tag{ii}\\ [B_{ij}, B_{jk}] = B^{-1}_{ik} \text{ and }[B_{ij},
B^{-1}_{jk}] = B_{ik}\tag{iii}\\ A_{ij} A^{-1}_{ji} A_{ij} =
A^{-1}_{ji} A_{ij} A^{-1}_{ji} \text{ and } B_{ij} B^{-1}_{ji} B_{ij} =
B^{-1}_{ji} B_{ij} B^{-1}_{ji}\tag{iv}
\end{gather}
\end{lem}
\begin{proof}
Let $\{e_i\}_{i=1}^n$ be the generators of 
$F_n$ and define automorphisms $A_{ij}$ and 
$B_{ij}$ by
\begin{align*} 
A_{ij}(e_k) &= 
\begin{cases}
e_i e_j,	&\text{if $i = k,$} \\ 
e_k,	&\text{ otherwise, and }\\
\end{cases}\\
B_{ij}(e_k) &= 
\begin{cases}
e_j e_i,	&\text{if $i = k,$} \\ 
e_k,	&\text{ otherwise.}
\end{cases}
\end{align*}

Then $\{A_{ij}, B_{ij}\}$ with $i \ne j,$ generate the index two subgroup of 
$\Aut( F_n)$ given by those automorphisms which induce 
on the abelianization $\Z^n$ of $F_n$ an automorphism of determinant
one (see, e.g.\ \cite{LS}).   
A straightforward but tedious computation shows that 
relations (i) -- (iv) are satisified.  
\end{proof}

We next find some fixed points.

\begin{lem}
\label{lem:fixed.generators}
Suppose $n>4$ and $\phi: \Aut( F_n) \to \Diff_+^2(\T^1)$ is a
homomorphism with $a_{ij} = \phi(A_{ij})$ and $b_{ij} = \phi(B_{ij})$
where $A_{ij}$ and $B_{ij}$ are the generators of Lemma
\ref{lem:generators}.  Then each of the diffeomorphisms $a_{ij}$ and
$b_{ij}$ has a fixed point.
\end{lem}
\begin{proof}
We fix $i, j$ and show $a_{ij}$ has a fixed point.
Since $n > 4$ there is an $a_{kl}$ with $\{i,j\} \cap \{k,l\} = \emptyset.$
Let $Z(a_{kl})$ denote the centralizer of $a_{kl}$, so $a_{ij} \in Z(a_{kl}).$
Also, since $n > 4$ there is $1 \le q \le n$ which is distinct from $i,j,k,l,$
so $a_{iq}, a_{qj} \in Z(a_{kl})$.

By Lemma \ref{lem:rho.centralizer} we know that the rotation number 
$\rho: Z(a_{kl}) \to \R$ is a homomorphism so 
$\rho(a_{ij})  = \rho( [a_{iq}, a_{qj}]) = 0.$  This implies that
$a_{ij}$ has a fixed point.
\end{proof}

We can now prove the main result of this section

\bigskip
\noindent
{\bf Theorem \ref{theorem:autogroups}.}
{\it 
For $n \ge 6$, any homomorphisms from $\Aut( F_n)$ or or $\Out( F_n)$ to
$\Diff_+^2(I)$  or $\Diff_+^2(\T^1)$ 
factors through $\Z/2\Z$.}
\bigskip

\begin{proof}
Let $H$ be the index two subgroup of $\Aut(F_n)$ from Lemma
\ref{lem:generators}.  Let $\phi$ be a homomorphism $H$ to
$\Diff_+^2(I)$ or $\Diff_+^2(\T^1)$ and let $a_{ij} = \phi(A_{ij})$ and
$b_{ij} = \phi(B_{ij})$ where $A_{ij}$ and $B_{ij}$ are the generators
from Lemma \ref{lem:generators}.  Then $\{a_{ij}, b_{ij}\}$ are
generators for $Image(\phi).$ We will show in fact that the subgroup
$G_a$ generated by $\{a_{ij}\}$ is trivial, as is the subgroup $G_b$
generated by $\{b_{ij}\}.$ The arguments are identical so we consider
only the set $\{a_{ij}\}.$ Since $n \ge 6$, we have by Lemma
\ref{lem:generators} that the commutavity graph for the generators
$\{a_{ij}\}$ of $G_a$ is connected.  Hence by Corollary
\ref{cor:comm.graph.T^1}, if these generators are fully supported then
$G_a$ is abelian.  But the relations $[a_{ij}, a^{-1}_{jk}] = a_{ik}$
then imply that $G_a$ is trivial.

Hence we may assume at least one of these generators is not fully
supported.  Let $J$ be an interval which is a maximal component of
fixed point sets for any of the $a_i$. More precisely, we choose $J$ so
that there is $a_{pq}$ with $J$ a component of $\Fix(a_{pq})$ and so
that there is no $a_{kl}$ such that $\Fix(a_{kl})$ properly contains
$J$.

We wish to show that $J$ is fixed pointwise by each $a_{ij}$.  In case
$\{p,q\} \cap \{i,j\} = \emptyset$ we note that at least the 
endpoints of $J$ are fixed by $a_{ij}$ because
$a_{pq}$ and $a_{ij}$ commute and the endpoints of $J$ are in
$\partial \Fix( a_{pq}).$ So Lemma \ref{lem:fix.component.I} or Lemma
\ref{lem:per.component} implies these endpoints are fixed by $a_{ij}.$

For the general case we first show $a_{ij}(J) \cap J \ne \emptyset.$
In fact there is no $a_{ij}$ with the property that $a_{ij}(J) \cap J
= \emptyset$ because if there were then the interval $J'$ defined to
be the smallest interval containing $\{a^n_{ij}(J)\}_{n \in \Z}$ is an
$a_{ij}$-invariant interval with no interior fixed points for
$a_{ij}$.  (In case $G_a \subset \Diff_+^2(\T^1)$ we need the fact that
$a_{ij}$ has a fixed point, which follows from Lemma
\ref{lem:fixed.generators}, in order to know this interval exists.)  But
since $n \ge 6$ there is some $a_{kl}$ which commutes with both $a_{ij}$
and $a_{pq}$, and it leaves both $J$ and $J'$ invariant.  Applying Kopell's
Lemma (Theorem \ref{thm:kopell}) to $a_{kl}$ and $a_{ij}$ we conclude that $J'
\subset \Fix(a_{kl})$, which contradicts the maximality of $J$.  Hence
we have shown that $a_{ij}(J) \cap J \ne \emptyset$ for all $i,j.$

Now by Lemma \ref{lem:generators} the relation $a_{pq} a^{-1}_{qp} a_{pq} =
a^{-1}_{qp} a_{pq} a^{-1}_{qp},$ holds so we may apply Lemma \ref{lem:aba} to
conclude that $J$ is fixed pointwise by $a_{qp}.$

We next consider the generator $a_{pk}.$ If $x_0$ is an endpoint of
$J$ then since $a_{qk}(J) \cap J \ne \emptyset$ at least one of
$a_{qk}(x_0)$ and $a^{-1}_{qk}(x_0)$ must be in $J$.  Hence the
relations $[a_{pq}, a_{qk}] = a^{-1}_{pk}$ and $[a_{pq}, a^{-1}_{qk}]
= a_{pk}$ imply that $a_{pk}(x_0) = x_0.$ This holds for the other
endpoint of $J$ as well, so $J$ is invariant under $a_{pk}$ for all
$k.$  The same argument shows $J$ is invariant under $a_{qk}$ for all
$k.$ But since $a_{pq}$ is the identity on $J$ we conclude
from $[a_{pq}, a_{qk}] = a^{-1}_{pk}$ that $a_{pk}$ is the identity
on $J$.  Similarly $a_{qk}$ is the identity on $J$.

Next the relation $a_{qk} a^{-1}_{kq} a_{qk} = a^{-1}_{kq} a_{qk}
a^{-1}_{kq},$ together with Lemma \ref{lem:aba} implies that $a_{kq}$
is the identity on $J$.  A similar argument gives the same result
for $a_{kp}$.  Finally, the relation $[a_{ip}, a^{-1}_{pj}] = a_{ij}$
implies that $a_{ij}$ is the identity on $J.$

Thus we have shown that any subgroup of $\Diff_+^2(I)$ or
$\Diff_+^2(\T^1)$ which is a homomorphic image of $H$, the index two subgroup
of $\Aut(F_n),$ has an
interval of global fixed points.  In the case of $\Diff_+^2(\T^1)$ we can
split at a global fixed point to get a subgroup of $\Diff_+^2(I)$ which
is a homomorphic image of $H$.  In the $I$ case we can
restrict the action to a subinterval on which the action has no global
fixed point.  But
our result then says that for the restricted action there is an interval
of global fixed points, which is a contradiction.  We conclude that
the subgroup of $G_a$ generated by $\{a_{ij}\}$ is trivial and
the same argument applies to the subgroup $G_b$ generated by
$\{b_{ij}\}.$  So $\phi(H)$ is trivial.  Since $\Aut(F_n)/H$ has
order two any homomorphism $\phi$ from $\Aut(F_n)$ to $\Diff_+^2(I)$
or $\Diff_+^2(\T^1)$ has an image whose
order is at most two.  The fact that $\Diff_+^2(I)$ has no elements
of finite order implies that in this case $\phi$ is trivial.

Since there is a natural homomorphism from $\Aut(F_n)$ onto the group
$\Out(F_n)$, the statement of Theorem \ref{theorem:autogroups} for 
$\Out(F_n)$ holds.
\end{proof}

\subsection{Finite index and other subgroups}
\label{section:aut2}
In this section we prove Theorem
\ref{theorem:finite:index}.  Recall that the Torelli group $T_{g,k}$ is
the subgroup of $\Mod(g,k)$ consisting of diffeomorphisms of the
$k$-punctured, genus $g$ surface $\Sigma_{g,k}$ which act trivially on
$H^1(\Sigma_{g,k};\Z)$.  Of course any finite index subgroup of
$\Mod(g,k)$ contains a finite index subgroup of $T_{g,k}$, so we need
only prove the thoerem for the latter group; that is, (3) implies (1).

So, let $H$ be any finite index subgroup of $T_{g,k}, g\geq 3, k\geq 0$
or of $\Aut(F_n)$ or $\Out(F_n)$ for $n\geq 6$.  
Suppose that $H$ acts faithfully as a group 
of $C^2$ diffeomorphisms of $I$ or $S^1$, so that no nontrivial element
of $H$ has an interval of fixed points (i.e. the action is fully
supported, in the terminology of \S\ref{subsection:fully:supported}).
In each case we claim that $H$ contains infinite order elements $a,b,c$
with the property that, for each $r\geq 2$, the elements $a^r,b^r$
generate an infinite, non-abelian group which commutes with $c^s$ for
all $s$.  

To prove the claim for $T_{g,k}$ one takes $a,b,c$ to be Dehn twists about
homologically trivial, simple closed curves with $a$ and $b$ 
having positive intersection number and $c$ disjoint from both.  
For $r$ sufficiently large the elements $a^r,b^r,c^r$ clearly lie in
$H$.  As powers of Dehn twists commute if and only if the
curves have zero intersection number, the claim is proved for $T_g,k$.  
A direct calculation using appropriate generators $A_{ij}$ as in Lemma
\ref{lem:generators} gives the same claim for $\Aut(F_n)$ and
$\Out(F_n)$ for $n>2$.

Theorem \ref{theorem:finite:index} in the case of $I$ now follows
immediately from Theorem \ref{thm:comm.graph.I} and the fact that the
group $<a^r,b^r,c^r>$ has connected commuting graph but is non-abelian.
For the case of $S^1$, note that if $c^r$ has a fixed point for any
$r>0$ so do $a^r,b^r$ by Lemma \ref{lem:fix.component.I}, so after
taking powers we are done by the same argument, this time applying 
Corollary \ref{cor:comm.graph.T^1}.  

If no positive power of $c$ has a
fixed point on $S^1$, then $c$ must have irrational rotation number, in
which case the centralizer of $c$ in $\Homeo(S^1)$ is abelian by 
\ref{prop:irrational.rotation}.  But $<a,b>$ lie in this centralizer and 
do not commute, a contradiction.

The proof of Theorem \ref{theorem:finite:index} for $\Aut(F_n)$ and
$\Out(F_n)$ is similar to the above, so we leave it to the reader.

\section{Baumslag-Solitar groups}

In this section we study actions of the groups 
$$\BS(m,n)=<a,b:ab^ma^{-1}=b^n>$$
For $n>m>1$ these groups are not residually finite, hence are not
linear (see, e.g.\  \cite{LS}).

\subsection{An analytic action of $BS(m,n)$ on $\R$}

In this subsection we construct an example of an analytic action 
of $BS(m,n),\ 1 <m <n$ on $\R.$   
It is straightforward to find diffeomorphisms of $\R$ which satisfy the
Baumslag-Solitar relation.

\begin{prop}[Diffeomorphisms satisfying the Baumslag-Solitar relation]
\label{prop:coverings}
Let $f_n: \T^1 \to \T^1$ be any degree $n$ covering and let $f_m: \T^1
\to \T^1$ be any degree $m$ covering.  Let $g_n$ and $g_m$ be lifts of
$f_n$ and $f_m$ respectively to the universal cover $\R$.  Then the
group of homeomorphisms of $\R$ generated by $g = g_m g_n^{-1}$ and the
covering translation $h(x) = x +1$ satisfies the relation $gh^m g^{-1} =
h^n.$
\end{prop}
\begin{proof}
We have 
\begin{align*}
	g_n^{-1}( h^{n}(x)) &= g_n^{-1}( x+n )\\
	&= g_n^{-1}(x) + 1\\
	&= h(g_n^{-1}(x)),
\end{align*}
and
\begin{align*}
	g_m( h(x)) &= g_m( x + 1)\\
	&= g_m( x ) + m\\
	&= h^m( g_m(x)).
\end{align*}
So
\begin{align*}
	g(h^n(x)) &= g_m( g_n^{-1}( h^n(x)))\\
	&= g_m(h(g_n^{-1}(x)))\\
	&= h^m(g_m (g_n^{-1}(x)))\\
	&= h^m(g(x)).
\end{align*}
\end{proof}

We will now show how to choose the diffeomorphisms constructed in 
Proposition \ref{prop:coverings} more carefully so that we will be able
to use a Schottky type argument to show that the diffeomorphisms will
satisfy no other relations.  Part of the difficulty will be that the
``Schottky sets'' will have to have infinitely many components.

We first construct, for each $j \ge 1,$ a $C^\infty$ diffeomorphism
$\Theta_j$ of $\R$ which is the lift of a map of degree $j$ on $\T^1.$
Let $a = 1/10.$ We choose a $C^\infty$ function $\Theta_j : (-a/2, 1 +
a/2) \to \R$ with the following properties:
\begin{align*}
	\Theta_j(x) &= \frac{a}{2} x \text{ for } x \in (-a/2, a/2),\\
	\Theta_j(x) &= j + \frac{a}{2} (x-1) \text{ for } \in (a, 1+ a/2),\\
	\Theta_j'(x) &> 0.
\end{align*}
Note that for $x \in (-a/2, a/2)$ we have $\Theta_j(x + 1) =
\Theta_j(x) + j.$ It follows that we can extend $\Theta_j$ to all of
$\R$ by the rule $\Theta_j(x + k) = \Theta_j(x) + jk.$ Thus the
extended $\Theta_j$ is the lift of a map of degree $j$ on $\T^1.$
This map is in fact a covering map since $\Theta_j$ is a diffeomorphism.

We shall be particularly interested in the intervals
$A = [-a, a],\ A^s = [-a, 0],\ A^u = [0, a],$ and
$C = [a, 1-a].$ Note that $[0,1] = A^u \cup C \cup (A^s + 1).$
Also by construction $\Theta_j(A^s) \subset A^s$ and 
$\Theta_j(C) \subset (A^s + j).$  More importantly we have the
following inclusions concerning the union of integer translates of
these sets.  For a set $X \subset \R$ and a fixed $i_0$, we will denote 
by $X+i_0$ the set $\{x + i_0 |\ x \in X\}$, by $X+\Z$ the set  
$\{x + i |\ x \in X,\ i \in \Z\}$, and by $X+m\Z$ the set 
$\{x + mi |\ x \in X,\ i \in m\Z\}$.

Suppose $k$ is not congruent to $0$ modulo $j.$ Then we have
\begin{align*}
	\Theta_j(A^s + \Z) &\subset A^s + j\Z,\\
	\Theta_j(C + \Z) &\subset A^s + j\Z, \\
	\Theta_j(A^u + {j\Z} + k) &\subset A^s + j\Z\\
	\Theta_j^{-1}(C + \Z) &\subset A^u + \Z.
\end{align*}

We now define two diffeomorphisms of $\R$.  Let $g_n(x) = \Theta_n(x)$
and let $g_m(x) = \Theta_m(x - 1/2) + 1/2.$ We define $B = A + 1/2,\ B^s
= A^s + 1/2,\ B^u = A^u + 1/2,$ and $C_2 = C + 1/2.$ Then from the
equations above we have
\begin{align}
	g_n(C + \Z) &\subset A^s + n\Z\notag\\
	g_n^{-1}(C + \Z) &\subset A^u + \Z\notag\\
	g_n^{-1}(A + k) &\subset A^u + {\Z} 
	\text{ if } k \notin n\Z \label{eqn1}\\
	g_m(C_2 + \Z) &\subset B^s + m\Z\notag\\
	g_m^{-1}(C_2 + \Z) &\subset B^u + \Z\label{eqn2}\\
	g_m^{-1}(B + k) &\subset B^u + {\Z} 
	\text{ if } k \notin m\Z \label{eqn3}
\end{align}

The diffeomorphisms $g_m$ and $g_n$ can be approximated by analytic
diffeomorphisms which still satisfy the equations above and are
still lifts of covering maps on $\T^1.$  More precisely, the function
$\phi(x) = g_n(x) - nx$ is a periodic function on $\R$ and may be
$C^1$ approximated by a periodic analytic function.  If we replace
$g_n$ by $\phi(x) + nx$ then the new analytic $g_n$ will have postive
derivative and be a lift of a degree $n$ map on $\T^1.$  If the 
approximation is sufficiently close, 
it will also still satisfy the equations
above.  In a similar fashion we may perturb $g_m$ to be analytic
while retaining its properties.

Note that $B + {\Z} \subset C + {\Z}$ and $A + {\Z} \subset C_2 + \Z$,
so we have the following key properties:
\begin{align}
	g_n(B + \Z) &\subset A^s + n\Z,\label{eqn4}\\
	g_n^{-1}(B + \Z) &\subset A^u + \Z \label{eqn5}\\
	g_m(A + \Z) &\subset B^s + m\Z \label{eqn6}\\ 
	g_m^{-1}(A + \Z) &\subset B^u + \Z\label{eqn7}
\end{align}

Define $g(x) = g_n (g_m^{-1}(x))$ and $h(x) = x + 1,$ and let
$G$ be the group generated by $g$ and $h$.  
Propostion \ref{prop:coverings} implies that
$g h^m g^{-1} = h^n$.  Hence there is a surjective 
homomorphism $\Phi:\BS(m,n)\rightarrow G$ sending
$b$ to $g$ and $a$ to $h$.  We need only show that $\Phi$ 
is injective.

\begin{lem}[normal forms]
\label{lem:canonical}
In the group $BS(m,n)$ with generators $b$ and $a$ satisfying
the relation $ba^mb^{-1} = a^n$ 
every nontrivial element can be written in the form
\[
a^{r_n}b^{e_n}a^{r_{n-1}}b^{e_{n-1}} \dots a^{r_1}b^{e_1}a^{r_0}
\]
where $e_i = \pm 1$ and $\ r_i \in \Z$ have the property that
whenever $e_i =1$ and $e_{i-1} = -1$ 
we have $r_i \notin m\Z$, and whenever $e_i =-1$ and $e_{i-1} =1$
we have $r_i \notin n\Z.$
\end{lem}

While Lemma \ref{lem:canonical} follows immediately from the normal form 
theorem for HNN extensions, we include a proof here since it is so
simple in this case.

\begin{proof}
If we put no restrictions on the integers $r_i$ and, in particular
allow them to be $0$, then it is trivially true that any element
can be put in  the form above.  If for any $i,\ e_i =1,\ e_{i-1} =-1$
and $r_i = mk \in m\Z$ then the fact that $ba^{mk}b^{-1} = a^{nk}$ allows
us to substitute and obtain another expression in the form above,
representing the same element of $\BS(m,n)$, 
but with fewer occurences of the terms $b$ and $b^{-1}$.  
Similarly if $\ e_i =-1,\ e_{i-1} =1$
and $r_i = nk \in n\Z$ then $b^{-1}a^{nk}b$ can be replaced with $a^{mk}$
further reducing the occurences of the terms $b$ and $b^{-1}$.  

These substitutions can be repeated at most a finite number of times
after which we have the desired form.
\end{proof}

We now show that $\BS(m,n)$ is a subgroup of $\Diff_+^\omega(\R)$.  
Note that $\Diff_+^\omega(\R)\subset
\Homeo_+(I)\subset \Homeo_+(\T^1)$, the first inclusion being induced by the 
one-point compactification of $\R$.  Hence the following implies 
Theorem \ref{theorem:bs:existence}.  

\begin{prop}
\label{proposition:bs:isomorphism}
The homomorphism $\Phi : BS(m,n) \to G$ defined by 
$\Phi(b) = g$ and $\Phi(a) = h$ is an isomorphism.
\end{prop}

\begin{proof}
By construction $\Phi$ is surjective; we now prove injectivity.  
To this end, consider an arbitrary nontrivial element written in the
form of Lemma \ref{lem:canonical}.  
After conjugating we may assume that $r_n = 0$ (replacing 
$r_0$ by $r_0 - r_n$).  So we must show that ($\Phi$ applied to) 
the element 
\[
\alpha = g^{e_n}h^{r_{n-1}}g^{e_{n-1}} \dots h^{r_1}g^{e_1}h^{r_0}
\]
acts nontrivially on $\R.$  If $n = 0$, i.e.\ if 
$\alpha = h^{r_0}$, then $\alpha$ clearly acts nontrivially if $r_0 \ne
0$, so we may assume $n \ge 1.$

Let $x$ be an element of the interior of $C \cap C_2.$ We will prove by
induction on $n$ that $\alpha( x) \in (A^s +n{\Z}) \cup (B^s + m{\Z}).$
This clearly implies $\alpha(x) \ne x.$

Let $n = 1.$  Then assuming $e_1 = 1$, 
we have $\alpha(x) = g( h^{r_0}(x)) = g_n(g_m^{-1}( h^{r_0}(x))).$  
But $x \in C_2 + \Z$ implies $h^{r_0}(x) = x + r_0 \in C_2 + \Z.$
So equation (\ref{eqn2}) implies $g_m^{-1}(x + r_0) \in B^u + \Z$, and
then equation (\ref{eqn4}) implies $g_n(g_m^{-1}(x + r_0)) \in A^s + n\Z.$
Thus $e_1 = 1$ implies $\alpha(x)\in A^s + n\Z.$
One shows similarly if $e_1= -1$ then $\alpha(x) = g( h^{r_0}(x)) \in 
B^s + m\Z.$

Now as induction hypothesis assume that 
\[
y = g^{e_k}h^{r_{k-1}}g^{e_{k-1}} \dots h^{r_1}g^{e_1}h^{r_0}(x)
\]
and that either
\begin{align*}
	e_k &= 1 \text{ and } y \in A^s + n\Z, \text{ or }\\
	e_k &= -1 \text{ and } y \in B^s + m\Z.
\end{align*}

We wish to establish the induction hypothesis for $k+1.$
There are four cases corresponding to the values of $\pm1$ for 
each of $e_k$ and $e_{k+1}.$  If $e_k = 1$ and
$e_{k+1}= 1$ then $y \in A^s + n\Z$ so $h^{r_i}(y) =
y + r_i \in A^s + \Z$ and equation (\ref{eqn7}) implies that
$g_m^{-1}( y + r_i) \in B^u + \Z$.  Hence by equation (\ref{eqn4}),
$g(h^{r_i}(y)) = g_n(g_m^{-1}(y + r_i)) \in A^s + n\Z$  as desired.
On the other hand if $e_k = 1$ and $e_{k+1}= -1$ then $r_k
\notin n\Z$.  Hence $y \in A^s + n\Z$ implies $h^{r_i}(y) = y + r_i
\in A^s + p$ with $p \notin n\Z$.  Consequently 
$y' = g_n^{-1}(y + r_i) \in A^u + \Z$ by
equation (\ref{eqn1}).  So equation (\ref{eqn6}) implies
$g^{-1}(h^{r_i}(y)) = g_m(y') \in B^s + m\Z.$  Thus we have verified
the induction hypothesis for $k+1$ in the case $e_k = 1.$

If $e_k = -1$ and $e_{k+1}= -1$ then $y \in B^s + m\Z$ so $h^{r_i}(y)
= y + r_i \in B^s + \Z$ and equation (\ref{eqn5}) implies that
$g_n^{-1}( y + r_i) \in A^u + \Z$.  Hence by equation (\ref{eqn6})
$g(h^{r_i}(y)) = g_m(g_n^{-1}(y + r_i)) \in B^s + m\Z$ as desired.
Finally for the case $e_k = -1$ and $e_{k+1}= 1$ we have $r_k \notin
m\Z$.  Hence $y \in B^s + m\Z$ implies $h^{r_i}(y) = y + r_i \in 
B^s + p$ with $p \notin m\Z$.  
Consequently $y' = g_m^{-1}(y + r_i) \in B^u + \Z$ 
by equation (\ref{eqn3}).  So equation (\ref{eqn4}) implies
$g(h^{r_i}(y)) = g_n(y') \in A^s + n\Z.$ Thus we have verified the
induction hypothesis for $k+1$ in the case $e_k = -1$ also.
\end{proof}

\noindent
{\bf Remark.} We showed in the proof of Proposition
\ref{proposition:bs:isomorphism} that any element of the form
\[
g^{e_n}h^{r_{n-1}}g^{e_{n-1}} \dots h^{r_1}g^{e_1}h^{r_0}
\]
acts nontrivially on $\R$ provided $e_i = \pm 1,\ r_i \in \Z,$ whenever
$e_i =1$ and $e_{i-1} = -1$ we have $r_i \notin m\Z$, and whenever $e_i
=-1$ and $e_{i-1} =1$ we have $r_i \notin n\Z$.  As a
corollary we obtain the following lemma which we will need later.

\begin{lem}
\label{lem:nontrivial}
In the group $BS(m,n)$ with generators $b$ and $a$ satisfying
the relation $ba^mb^{-1} = a^n$ 
every element of the form
\[
b^{e_n}a^{r_{n-1}}b^{e_{n-1}} \dots a^{r_1}b^{e_1}a^{r_0},
\]
is nontrivial provided
$e_i = \pm 1,\ r_i \in \Z,$ and whenever $e_i =1 ,\ e_{i-1} = -1$
we have $r_i \notin m\Z$ and whenever $e_i =-1,\ e_{i-1} =1$
we have $r_i \notin n\Z.$
\end{lem}

Of course this lemma can also be obtained from the normal form theorem
for HNN extensions.

\subsection{General properties of $BS(m,n)$ actions}

In this short subsection we note two simple implications of the
Baumslag-Solitar relation.

\begin{lem}
\label{lem:BSperiodic}
If $BS(m,n)$ acts by homeomorphisms on $\T^1$ with generators $g,h$ satisfying
$gh^mg^{-1} = h^n$ then $h$ has a periodic point whose period
is a divisor of  $|n-m|.$
\end{lem}
\begin{proof} Consider the group of all lifts to $\R$ of elements of $BS(m,n)$ acting on $\T^1$.
Choose lifts $H$ and $G$ of $h$ and $g$ respectively and let $T(x) =
x+1$ denote the covering translation.  We have
\[
	 G H^m G^{-1} = H^n T^p
\]
for some integer $p.$ Each lift of an element of $BS(m,n)$ has a well
defined translation number in $\R$ and these are topological conjugacy
invariants. Hence
\[
 m \tau( H) = \tau( H^m) = \tau( G H^m G^{-1}) = \tau( H^n T^p) = p + n \tau(H).
\]
Solving we conclude
\[
	  \tau( H) = \frac{p}{m-n}.
\]
Since $H$ has a rational rotation number $h$ has a point whose period is a divisor of  $|n-m|.$
\end{proof}

\begin{lem}
Suppose $g$ and $h$ are homeomorphisms of $\R$ and
satisfy $gh^mg^{-1} = h^n.$   Then if $h$ is fixed point free, $g$
has a fixed point.
\end{lem}

\begin{proof}
Assume, without loss of generality, that $h(x) = x + 1$ and $n > m$
and $g(0) \in [p, p+1].$ Then 
\[
h^{nk}( g(0)) = (g(0) + nk) \in [p +nk, p+ nk + 1].
\]
But
\begin{align*}
	h^n g  &= g h^m, \text{ so }\notag\\
	h^{nk} g  &= g h^{mk}, \text{ and }\notag\\
	h^{nk}( g(0)) &= g ( h^{mk}(0)) = g(mk), \text{ so}\notag \\
	g(mk) &\in [p +nk, p+ nk + 1].\notag
\end{align*}
For $k$ sufficiently large this implies $g(mk) > mk$ and for $k$ sufficiently negative
that $g(mk) < mk.$  The intermediate value theorem implies $g$ has a fixed point.
\end{proof}

\subsection{$C^2$ actions of $BS(m,n)$}
In contrast to the analytic actions of $\BS(m,n)$ on $\R$, we will show
that there are no such actions, even $C^2$ actions, on either $I$ or
$\T^1$. 

\begin{lem}
\label{lem:bs.commute}
Suppose $g$ and $h$ are orientation-preserving $C^2$ diffeomorphisms
of $I$ satisfying $g h^m g^{-1} = h^n.$ Then $h$ and $ghg^{-1}$
commute.
\end{lem}
\begin{proof}
If $x \in \Fix(h)$ then $g(x) = g(h^m(x)) = h^n(g(x))$ so $g(x)$ is a periodic
point for $h.$  But the only periodic points of $h$ are fixed.  We conclude
that $g(\Fix(h))= \Fix(h).$

Let $(a,b)$ be any component of the complement of $\Fix(h)$ in $[0,1].$
Since $g(\Fix(h)) = \Fix(h)$ we
have that $ghg^{-1}([a,b]) = [a,b].$ 

Let $H$ be the centralizer of $h^n = gh^mg^{-1}$.  Note that $h$ and
$ghg^{-1}$ are both in $H.$ Let $f$ be any element of $H$.  Then if
$f$ has a fixed point in $(a,b)$, since it commutes with $h^n$ we know
by Kopell's Lemma (Theorem \ref{thm:kopell}) 
that $f = id.$ In other words $H$ acts
freely (though perhaps not effectively) on $(a,b).$ Hence the
restrictions of any elements of $H$ to $[a,b]$ commute by H\"older's Theorem.
In particular the restrictions of $h$ and $ghg^{-1}$ to $[a,b]$ commute.
But $(a,b)$ was an arbitrary component of the complement of $\Fix(h)$.
Since the restrictions of $h$ and $ghg^{-1}$ to $\Fix(h)$ are both the
identity we conclude that $h$ and $ghg^{-1}$ commute on all of $I$.
\end{proof}

As a corollary we have:

\bigskip
\noindent
{\bf Theorem \ref{theorem:bs:noact} for $I$: }{\it 
If $m$ and $n$ are greater than $1$ then no subgroup of $\Diff_+^2(I)$
is isomorphic to $BS(m,n).$}
\bigskip

\begin{proof}
If there is such a subgroup, it has generators $a$ and $b$ satisfying
$a^n = ba^mb^{-1}$.  Lemma \ref{lem:bs.commute} asserts that
$a$ and $c = bab^{-1}$ commute.  But the commutator $[c,a]$ is
$ba^{-1}b^{-1}abab^{-1}a$ and this element is nontrivial by
Lemma \ref{lem:nontrivial}.  This contradicts the existence of
such a subgroup.
\end{proof}

Showing that $\BS(m,n)$ has no $C^2$ action on the circle is a bit
harder.

\bigskip
\noindent
{\bf Theorem \ref{theorem:bs:noact} for $\T^1$: }{\it 
If $m$ is not a divisor of $n$ then no subgroup of $\Diff_+^2(\T^1)$
is isomorphic to $BS(m,n).$}
\bigskip

\begin{proof}
Suppose there is such a subgroup, and we have generators $g,h$
satisfying $gh^mg^{-1} = h^n$.  By Lemma \ref{lem:BSperiodic}, $h$ has a
periodic point.  Let $Per(h)$ be the closed set of periodic points of
$h$ (all of which have the same period, say $p$).  If $x \in \Per(h)$
then $g(x) = gh^{mp} = h^{np}g(x)$ which shows $g(\Per(h)) = \Per(h).$

The omega limit set $\omega(x,g)$ of a point $x$ under $g$ is equal to a
subset of the set of periodic points of $g$ if $g$ has periodic points,
and is independent of $x$ if $g$ has no periodic points.  Since $g$ is
$C^2$ we know by Denjoy's theorem that $\omega(g, x)$ is either all of
$\T^1$ or is a subset of $\Per(g)$ (see \cite{dMvS}).  In the case at hand 
$\omega(g,x)$ cannot be all of
$\T^1$ since $x \in \Per(h)$ implies $\omega(x,g) \subset \Per(h)$ and
this would imply $h^p = id.$ Hence there exists $x_0 \in \Per(g) \cap
\Per(h).$ Then $x_0 \in \Fix(h^p) \cap \Fix( g^q)$ where $q$ is the
period of points in $\Per(g).$

The fact that $x_0$ has period $p$ for $h$ implies that its rotation
number under $h$ (well defined as an element of $\R/\Z)$ is $k/p + \Z$
for some $k$ which is relatively prime to $p$.  Using the fact that
$\rho(h^m) = \rho( g h^m g^{-1})$,  we conclude that 
$$m(k/p) + \Z = \rho( h^m)=\rho( g h^m g^{-1}) = \rho( h^n) = n(k/p) +
\Z$$
Consequently $mk \equiv nk\mod p$ and since $k$ and $p$ are relatively
prime $m \equiv n\mod p.$ The fact that $1 < m < n$ and $n-m = rp$ for
some $r \in \Z$ tells us that $p \notin m \Z$ because if it were then
$n$ would be a multiple of $m$.  Similarly $p \notin n \Z.$

Let $G$ be the group generated by $h_0 = h^p$ and $g_0=g^q$.  Now $G$ has a
global fixed point $x_0$, so we can split $\T^1$ at $x_0$ to obtain a $C^2$
action of $G$ on $I$. Note that $g_0h_0^{m^q}g_0^{-1} = g^q
h^{pm^q}g^{-q} = h^{pm^q} = h_0^{n^q}$.

We can apply Lemma \ref{lem:bs.commute} to conclude that $h_0$ and
$g_0h_0g_0^{-1}$ commute.  Thus $h^p$ and $c = g^qh^pg^{-q}$ commute.
But their commutator
\[
[c,h^p] =  g^qh^{-p}g^{-q}h^{-p}g^qh^pg^{-q}h^p.
\]
Since $p \notin m \Z$ and $p \notin n \Z$ this element is nontrivial
by Lemma \ref{lem:nontrivial}.  This contradicts the existence of
such a subgroup. 
\end{proof}

\subsection{Local rigidity of the standard action of $BS(1,n)$}

We consider the action of $BS(1,n)$ on $\R$ given by
elements of the form $x \to n^{k} x + b$ where $k \in \Z$ and
$b \in \Z[1/n].$  It is generated by the two elements 
$g_0(x) = nx$ and $h_0(x) = x + 1.$ 

The following result is essentially the theorem of Shub from
\cite{Sh}.   We give the proof since framing his result in our
context is nearly as long.

\begin{thm}
\label{theorem:local_rigid}
There are neighborhoods of $g_0$ and $h_0$ in the uniform 
$C^1$ topology such that whenever $g$ and $h$ are chosen from these
respective neighborhoods, and when the correspondences $g_0\rightarrow
g$ and $h_0\rightarrow h$ induce an isomorphism of $BS(1,n)$ to the group 
generated by $g$ and $h$, then the perturbed action is topologically
conjugate to the original action.
\end{thm}

\begin{proof}
Since $g$ is uniformly $C^1$ close to $g_0(x) = nx$, its inverse is a
contraction of $\R$ and has a unique fixed point, which after a change
of coordinates we may assume is $0$.  Since $h$ is uniformly close to
$h_0(x) = x+1$ it is fixed point free and after a further change of
co-ordinates we may assume $h^n(0) = n.$ Consider the space $H$ of
$C^0$ maps $\phi: \R \to \R$ such that $h_0\phi =\phi h.$ and $h(0) =
0.$ Any $\phi \in H$ satisfies $\phi(n) = n$ for all $n \in \Z.$
Moreover, such a $\phi$ is completely determined by its values on the
interval $[0,1]$, since $\phi(x + n) = \phi(x) +n.$ Indeed any $C^0$
map $\phi: [0,1] \to \R$ is the restriction of some element of $H$.
Clearly, a sequence in $H$ will converge if and only if it converges
on $[0, 1].$ Hence $H$ can be considered a closed subset of the
complete metric space of continuous functions on $[0, 1]$ with the
$C^0$ {\it sup} norm.

Consider the map  $F$ on $H$ given by $G(\phi) = g_0^{-1} \phi g$.  
Then $G(\phi) h = g_0^{-1} \phi g h = g_0^{-1} \phi h^n g
= g_0^{-1} h_0^n \phi  g = h_0 g_0^{-1} \phi  g = h_0 F(\phi).$
It follows $F: H \to H.$  Moreover in the $C^0$ {\it sup} norm this map
is easily seen to be a contraction.  It follows that there is a unique
fixed point $\phi_0$ for $F$ in $H$.

The map $\phi_0: \R \to \R$ is continuous and satisfies
\begin{align*}
	h_0^n \phi_0 &= \phi_0 h^n, \text{ and }\\
	g_0^n \phi_0 &= \phi_0 g^n,
\end{align*}
for all $n \in \Z$.  The first of these equations implies $\phi_0$
has image which is unbounded above and below, so $\phi_0$ is surjective.
This equation also implies that there is a uniform bound on 
the size of $\phi_0^{-1}(x)$ for $x\in \R.$  Indeed if $M$ is an integer
and $M > sup | \phi_0(0)- \phi_0(x)|$ for $x \in [0, 1]$ then 
$y > M+1$ or $y < -M$ implies $\phi_0(y) \ne \phi_0(x).$
So the sets $\phi_0^{-1}(z)$ have diameters with an upper bound independent
of $z \in \R.$  In particular, $\phi_0$ is proper.

We can now see that $\phi_0$ is injective.   Since 
$\phi_0(y) = \phi_0(x)$  implies $g_0^n(\phi_0(y)) = g_0^n(\phi_0(x))$
which in turn implies $\phi_0(g^n(y)) = \phi_0(g^n(x))$ and $g$ is
is a uniform expansion, it follows that if $\phi_0$ fails to be injective
then the sets $\phi_0^{-1}(z)$ do not have a diameter with an upper
bound independent of $z.$  This contradiction implies $\phi_0$ must
be injective.  The fact that $\phi_0$ is proper implies $\phi_0^{-1}$
is continuous.  Hence $\phi_0$ is a topological conjugacy from the
standard affine action of $BS(1,n)$ on $\R$ to the action of the
group generated by $g$ and $h.$
\end{proof}

In contrast to local rigidity of the standard action, Hirsch \cite{H} has found
real-analytic actions of $\BS(1,n)$ on the line which are not 
topologically conjugate to the standard action.

\begin{thm}[Hirsch \cite{H}]
\label{theorem:no_global_rigid}
There is an analytic action of $BS(1,n)$ on $\R$ which is not topologically
conjugate to the standard affine action.
\end{thm}

The construction is similar to our construction of an analytic $BS(m,n)$
action on $\R$.  The relation $ghg^{-1} = h^n$ is satisfied by analytic
diffeomorphisms $g$ and $h$ and $g$ has a fixed attracting fixed point.
Hence $g$ cannot be topologically conjugate to the affine function
$g_0(x) = nx.$

\bigskip
\noindent
Benson Farb:\\
Dept. of Mathematics, University of Chicago\\
5734 University Ave.\\
Chicago, Il 60637\\
E-mail: farb@math.uchicago.edu
\medskip

\noindent
John Franks:\\
Dept. of Mathematics, Northwestern University\\
Evanston, IL 60208\\
E-mail: john@math.northwestern.edu

\end{document}